\documentclass [12pt] {report}
\usepackage[dvips]{graphicx}
\usepackage{epsfig, latexsym,graphicx}
\usepackage{amssymb}
\newcommand {\D}[2] {\displaystyle\frac{\partial{#1}}{\partial{#2}}}
\newcommand {\DD}[2] {\displaystyle\frac{\partial^2{#1}}{\partial{#2^2}}}

\newcommand {\Dd}[3] {\displaystyle\frac{\partial^2{#1}}{\partial{#2}\partial{#3}}}

\newcommand {\ga} {\gamma}
\newcommand {\la} {\lambda}

\newcommand {\de} {\delta}
\newcommand {\prtl} {\partial}
\newcommand {\fr} {\displaystyle\frac}

\newcommand {\be} {\begin{equation}}
\newcommand {\ee} {\end{equation}}
\newcommand {\ba} {\begin{array}}
\newcommand {\ea} {\end{array}}
\newcommand {\bp} {\begin{picture}}
\newcommand {\ep} {\end{picture}}
\newcommand {\bc} {\begin{center}}
\newcommand {\ec} {\end{center}}
\newcommand {\bt} {\begin{tabular}}
\newcommand {\et} {\end{tabular}}
\newcommand {\lf} {\left}
\newcommand {\rg} {\right}

\newcommand {\cF} {{\cal F}}

\newcommand {\cR} {{\cal R}}

\newcommand {\ses} {\medskip}

\newcommand {\e} {\mathop{\rm e}\nolimits}

\newcommand {\bibit} {\bibitem}
\newcommand {\nin} {\noindent}

\usepackage{psfrag}

\usepackage{epsfig, latexsym,graphicx}
 \usepackage{amsmath}

\def\2#1#2#3{{#1}_{#2}\hspace{0pt}^{#3}}
\def\3#1#2#3#4{{#1}_{#2}\hspace{0pt}^{#3}\hspace{0pt}_{#4}}
\newcounter{sctn}
\def\sec#1.#2\par{\setcounter{sctn}{#1}\setcounter{equation}{0}
                  \noindent{\bf\boldmath#1.#2}\bigskip\par}

\begin {document}

\begin {titlepage}

\vspace{0.1in}

\begin{center}

{\Large \bf  Finsleroid-regular    space  developed.    Berwald case}

\end{center}

\vspace{0.3in}

\begin{center}

\vspace{.15in} {\large G.S. Asanov\\} \vspace{.25in}
{\it Division of Theoretical Physics, Moscow State University\\
119992 Moscow, Russia\\
{\rm (}e-mail: asanov@newmail.ru{\rm )}} \vspace{.05in}

\end{center}

\begin{abstract}

\ses

The Finsleroid--Finsler space becomes regular  when
the norm $||b||=c$ of the input 1-form  $b$ is taken to be  an arbitrary
positive scalar $c(x) <  1$.  By performing required direct evaluations,
  the respective spray coefficients have  been obtained in a simple and transparent form.
The  adequate    continuation  into  the  regular  pseudo-Finsleroid  domain  has  been  indicated.
The Finsleroid-regular Berwald space is found under the assumptions that the Finsleroid charge is a constant and
the 1-form $b$ is  parallel.

\ses

\nin{\it Keywords:} Finsler  metrics, spray coefficients, curvature tensors.

\end{abstract}

\end{titlepage}

\vskip 1cm

\ses

\ses

\setcounter{sctn}{1} \setcounter{equation}{0}

\nin
  {\bf 1. Description  of new conclusions}

\ses

\ses

In the Finsleroid-Finsler space $\cF\cF^{PD}_g $
(as well in its relativistic counterpart space $\cF\cF^{SR}_g $)
 constructed and studied in the previous papers [1-4]
the consideration was referred everywhere to
 the case $||b||=1$, that is,  the  vector field
 $b_i(x)$ involved in the 1-form $b=b_i(x)y^i$ was assumed to be of  unit length.
In the present paper,
we expand the restrictive case to get the possibility $||b||<1$.

Namely, we shall  deal with the Finsler space notion specified by the condition
that the Finslerian metric function
$K(x,y)$ be of the functional dependence
\be
K(x,y) =\Phi \Bigl(g(x), b_i(x), a_{ij}(x),y\Bigr)
\ee
subject to the conditions
\be
a^{ij}(x)b_i(x)b_j(x)=c^2(x)
\ee
and
\be
0<c(x)<1.
\ee
In (1.1),   $ a_{ij}(x)$
is a Riemannian metric tensor and
$g(x)$ plays the role of  the Finsleroid charge.
The explicit form of the $\Phi$ is specified by the representations written out explicitly
in Appendix A   in the positive-definite case,  and
in Appendix B   in the relativistic  case, respectively
 we obtain  the spaces ${\mathbf\cF\cR^{PD}_{g;c} } $  and  ${\mathbf\cF\cR^{SR}_{g;c} } $.
They fulfill the correspondence
\be
\cF\cR^{PD}_{g;c=1}    =   \cF\cF^{PD}_g
\ee
and
\be
\cF\cR^{SR}_{g;c=1}    =   \cF\cF^{SR}_g.
\ee

The Riemannian  squared metric  function
\be
S^2=b^2+q^2
\ee
 underlines the  positive-definite Finsler space  ${\mathbf\cF\cR^{PD}_{g;c} } $
 under study, and the pseudo-Riemannian
$  S^2=b^2-q^2  $ is to relate to the space ${\mathbf\cF\cR^{SR}_{g;c} } $.

\ses


This  scalar    $c(x)$ proves to play the role of the {\it regularization factor.}
Indeed, in the space ${\mathbf\cF\cR^{PD}_{g;c} } $
the metric function $K$   is  constructed in accordance with the formulas (A.19)-(A.23)  which involve
 the square-root variable
$q(x,y)=\sqrt{S^2-b^2}$ (see (1.6)).
Differentiating various  tensors of  the space ${\mathbf\cF\cR^{PD}_{g;c} } $ gives rise,
therefore,
 to appearance of degrees of the fraction $1/q$. Because of the inequality (A.5), the fractions $1/q$
 do not  produce any singularities as far as $ c<1$. If, however, $c=1$, we have $q(x,y)=0$ when the vector
 $y\in T_xM$ is a factor of the vector $b^i(x)$, so that in the space $\cF\cF^{PD}_g$
 the singularities can appear  when the vector $y$ belongs to the Finsleroid axis.

In contrast, in the space $\cF\cR^{PD}_{g;c} $ which uses $c\ne1$
the variable $q$ does not vanish anywhere
over all the slit tangent bundle $TM\setminus 0$,
and we obtain
the function $\Phi$ which is smooth  of the  class $C^{\infty}$.

We use

\ses

\ses

{REGULARITY   DEFINITION.}   The Finsler space under study is {\it regular} in the following sense:
over all the slit tangent bundle $TM\setminus 0$,
the Finsler metric function $K(x,y)$ of the space is smooth of the class $C^{\infty}$
and the entailed Finsler metric tensor $g_{ij}(x,y)$ is positive-definite: $\det( g_{ij}) >0$.

\ses

\ses

Attentive direct calculations
of
the induced spray coefficients
$ G^i =\ga^i{}_{nm}y^ny^m$,
where
$\ga^i{}_{nm}$ denote the associated Finslerian Christoffel symbols,
can be used to arrive at the following result.


\ses
\ses

THEOREM 1.  {\it In the Finsleroid-regular  space
 ${\mathbf\cF\cR^{PD}_{g;c} } $
the spray coefficients $G^i$
can  explicitly be written in the form}
\be
G^i=
\fr  g{\nu}
\Bigl(
y^jy^h\nabla_jb_h
+gqb^jf_j\Bigr)
v^i
-gqf^i
 +E^i
 +a^i{}_{nm}y^ny^m.
\ee

\ses

\ses

We use  the  notation
\be
v^i=y^i-bb^i
\ee
and
\be
f_j=f_{jn}y^n,\quad
f^i=f^i{}_ny^n,\quad
f^i{}_n=a^{ik}f_{kn}, \quad
f_{mn}=
\nabla_mb_n-\nabla_nb_m
\equiv \D{ b_n}{x^m}-\D {b_m}{x^n},
\ee
where
  $\nabla$ means the covariant derivative in terms of the associated Riemannian space
${R}_N=(M,{S})$  (see (A.83));
  $ a^i{}_{nm} $ stands for
the  Riemannian  Christoffel symbols (A.84) constructed from the input Riemannian metric tensor $a_{ij}(x)$;
the coefficients $E^i$ involving the gradients
 $g_h=\partial g/\partial x^h$
 of the Finsleroid charge
 can be taken as
\be
E^i = M (yg)y^i
+  K \fr{2 b^2 w^2}{gB}(yg) X A^i
-\fr12 MK^2g_hg^{ih},
\ee
where $(yg)=g_hy^h$;
$X$ is the function given in (A.67);  $w$ is the variable (A.38);
the function $M$ entered (1.10) is defined by $\partial K/\partial g= (1/2)MK$.

\ses

\ses

The difference
$
G^i-E^i-a^i{}_{nm}y^ny^m
$
involves  the  crucial terms linear in the covariant derivative
 $ \nabla_jb_h.
 $

Also, the following theorem is valid.

\ses
\ses

THEOREM 2.  {\it In the Finsleroid-regular  space
 ${\mathbf\cF\cR^{PD}_{g;c} } $
the equality
 \be
 G^i=  a^i{}_{nm}y^ny^m
\ee
 holds if and only if
\be
g=const ~ ~ { \rm and}  ~  ~
\nabla_mb_n=0.
\ee
}

\ses

\ses

When the equality (1.11) holds, one says that the Finsler space is the Berwald space (see [5]).
The above theorem yields a simple and attractive example of the regular Berwald space.

The sufficiency of the conditions (1.12) is obvious from the representation (1.7).
To verify the necessity, it is worth noting  that in the Berwald case the covariant derivative of the Cartan tensor
$A_{ijk}$ vanishes identically (see [5,6]), in which case the representations  (A.71) and (A.72) just entail $g=const$,
which in turn yields $E^i=0$.
With this observation, it is easy to see that the representation (1.7) reduces to (1.11) if only  $\nabla_mb_n=0$,
as far as the Finsleroid charge $g$ is kept differing from zero (the choice $g=0$ would reduce the Finsleroid-Finsler
 space to a Riemannian space).

\ses

\ses\ses

\setcounter{equation}{0}

{
\nin \bf Appendix A:   Distinguished
 ${\mathbf\cF\cR^{PD}_{g;c} } $-notions
 }

\ses\ses

Let $M$ be an $N$-dimensional
$C^{\infty}$
differentiable  manifold, $ T_xM$ denote the tangent space to $M$ at a point $x\in M$,
and $y\in T_xM\backslash 0$  mean tangent vectors.
Suppose we are given on $M$ a Riemannian metric ${\cal S}=S(x,y)$.
 Denote by
$\cR_N=(M,{\cal S})$
the obtained $N$-dimensional Riemannian space.
Let us also assume that the manifold $M$ admits a non--vanishing 1-form $b= b(x,y)$,
denote by
\be
c=|| b||_{\text{Riemannian}}
\ee
the respective Riemannian norm value.
Assuming
\be
 0 < c < 1,
 \ee
we get
\be
S^2-b^2>0
\ee
and may conveniently use the variable
\be
q:=\sqrt{S^2-b^2}.
\ee
Obviously, the inequality
\be
q^2 \ge \fr{1-c^2}{c^2}\,b^2
\ee
is valid.

With respect to  natural local coordinates in the space
$\cR_N$
we have the local representations
\be
a^{ij}b_ib_j=c^2
\ee
and
\be
 b=b_i(x)y^i,  \qquad S= \sqrt{a_{ij}(x)y^iy^j}.
\ee
The reciprocity  $a^{in}a_{nj}=\de^i{}_j$ is assumed, where $\de^i{}_j$ stands for the Kronecker symbol.
The covariant index of the vector $b_i$  will be raised by means of the Riemannian rule
$ b^i=a^{ij}b_j,$ which inverse reads $ b_i=a_{ij}b^j.$
We
also  introduce the tensor
\be
r_{ij}(x)~:=a_{ij}(x)-b_i(x)b_j(x)
\ee
to have the representation
\be
q=\sqrt{r_{ij}(x)y^iy^j}.
\ee
The equalities
\be
r_{ij}b^j=(1-c^2)b_i, \qquad
  r_{in}r^{nj}=r^j{}_i-(1-c^2)b^jb_i
 \ee
hold.

We introduce on the background manifold $M$
a scalar field $g=g(x)$
 subject to ranging
\be
-2<g(x)<2,
\ee
and apply  the convenient notation
\be
h(x)=\sqrt{1-\fr14(g(x))^2}, \qquad
G(x)=\fr{g(x)}{h(x)}
\ee
(compare with (2.10) in [3]).
The Finsleroid-regular space is underlined by
the {\it  characteristic  quadratic form}
\be
B(x,y) :=
b^2+gqb+q^2
\equiv
\fr 12
\Bigl[(b+g_+q)^2+(b+g_-q)^2\Bigr]
\ee
(cf.  (2.11) in [3]),
where $ g_+=\fr12g+h$ and $ g_-=\fr12g-h$,
 and the  discriminant $D_{\{B\}}$ of the quadratic form $B$
is negative:
\be
D_{\{B\}}   =    -4h^2<0.
\ee
Therefore, the {\it quadratic form $B$ is  positively definite.}
In the limit $g\to 0$,
the definition (A.13) degenerates to the
 quadratic form  of the input Riemannian metric tensor:
\be
 B|_{_{g=0}}=b^2+q^2 \equiv S^2.
\ee
Also,
\be
\eta B|_{_{y^i=b^i}}=c^2,
\ee
where
\be
\eta=\fr 1{1+gc\sqrt{1-c^2}}.
\ee
On the definition range (A.11) of the $g$,  we have
\be
\eta>0.
\ee

Under these conditions,
we introduce the following definition.

\ses

\ses

 { DEFINITION}. The scalar function $K(x,y)$ given by the formulas
\be
K(x,y)=
\sqrt{B(x,y)}\,J(x,y)
\ee
and
\be
J(x,y)=\e^{-\frac12G(x)f(x,y)},
\ee
where
\be
f=
-\arctan
 \fr G2
+\arctan\fr{L}{hb},
\qquad  {\rm if}  \quad b\ge 0,
\ee
and
\be
f= \pi-\arctan
\fr G2
+\arctan\fr{L}{hb},
\qquad  {\rm if}
 \quad b\le 0,
\ee
 with
 \be
 L =q+\fr g2b,
\ee
\ses\\
is called
the {\it  Finsleroid-regular   metric function}.

\ses

\ses

The function (A.23) obeys   the identity
\be
L^2+h^2b^2=B.
\ee

\ses

\ses

 { DEFINITION}.  The arisen  space
\be
\cF\cR^{PD}_{g;c} :=\{\cR_{N};\,b_i(x);\,g(x);\,K(x,y)\}
\ee
is called the
 {\it Finsleroid-regular  space}.

\ses

\ses

 { DEFINITION}. The space $\cR_N$ entering the above definition is called the {\it associated Riemannian space}.

\ses

\ses

The associated Riemannian metric tensor $a_{ij}$ has the meaning
\be
 a_{ij}=g_{ij}\bigl|_{g=0}\bigr. .
\ee

\ses

\ses

 { DEFINITION}.\, Within  any tangent space $T_xM$, the Finsleroid-regular metric function $K(x,y)$
 produces the {\it regular Finsleroid}
 \be
 \cF\cR^{PD}_{g;c\, \{x\}}:=\{y\in   \cF\cR^{PD}_{g;c\, \{x\}}: y\in T_xM , K(x,y)\le 1\}.
  \ee

\ses

 \ses

 { DEFINITION}.\, The {\it regular Finsleroid Indicatrix}
 $ I\cR^{PD}_{g;c\, \{x\}} \subset T_xM$ is the boundary of the regular Finsleroid, that is,
 \be
 I\cR^{PD}_{g;c\, \{x\}} :=\{y\in I\cR^{PD}_{g;c\, \{x\}} : y\in T_xM, K(x,y)=1\}.
  \ee

\ses

 \ses

 {\large  Definition}. The scalar $g(x)$ is called
the {\it Finsleroid charge}.
The 1-form $b=b_i(x)y^i$ is called the  {\it Finsleroid--axis}  1-{\it form}.

\ses

\ses

We shall meet the function
\be
\nu:=q+(1-c^2)gb
\ee
for which
\be
\nu >0 \quad \text{when} \quad |g|<2.
\ee
Indeed, if $gb>0$, then  the right-hand part of (A.29) is positive.
 When $gb<0$, we may note that
 at any fixed $c$ and $b$ the minimal value of $q$ equals
$\sqrt{1-c^2}|b|/c$ (see (A.5)), arriving again at (A.30).

The identities
\be
\fr{c^2S^2-b^2}{q\nu}=1-(1-c^2)\fr B{q\nu},
\qquad
gb(c^2S^2-b^2)=qB-\nu S^2
\ee
are valid.

In evaluations it is convenient to use the variables
\be
u_i~:=a_{ij}y^j,
\qquad
v^i~:=y^i-bb^i, \qquad v_m~:=u_m-bb_m=r_{mn}y^n\equiv a_{mn}v^n.
\ee
We have
\be
r_{ij}=\D{v_i}{y^j},
\ee
\ses
\be
u_iv^i=v_iy^i=q^2,
 \qquad
v_ib^i=v^ib_i=(1-c^2)b,
\ee
\ses
\be
 r_{in}v^n=v_i-(1-c^2)bb_i,  \qquad
v_kv^k=q^2-(1-c^2)b^2,
\ee
and
\be
\D b{y^i}=b_i, \qquad \D q{y^i}=\fr{v_i}q.
\ee

Under these conditions,
 we are to  explicitly extract from the function $K$ the  distinguished Finslerian tensors,
 and first of all
the covariant tangent vector $\hat y=\{y_i\}$,
the  Finslerian metric tensor $\{g_{ij}\}$
together with the contravariant tensor $\{g^{ij}\}$ defined by the reciprocity conditions
$g_{ij}g^{jk}=\de^k_i$, and the  angular metric tensor
$\{h_{ij}\}$, by making  use of the following conventional  Finslerian  rules in succession:
\be
y_i :=\fr12\D{K^2}{y^i}, \qquad
g_{ij} :
=
\fr12\,
\fr{\prtl^2K^2}{\prtl y^i\prtl y^j}
=\fr{\prtl y_i}{\prtl y^j}, \qquad
h_{ij} := g_{ij}-y_iy_j\fr1{K^2}.
\ee

To this end it is convenient to use the variable
\be
w=\fr qb,
\ee
obtaining
\be
\D{w}{y^i}=\fr{z_i}{b^2q}, \qquad z_i=bv_i-q^2b_i
\equiv bu_i-S^2b_i,
\ee
and
$$
y^iz_i=0,\qquad
b^iz_i=b^2-c^2S^2,
$$
\ses
$$
 a^{ij}z_iz_j=S^2(c^2S^2-b^2).
$$
We also introduce  the $\eta$--{\it tensors}  given by
\be
\eta_{ij}~:=r_{ij}-\fr1{q^2}v_iv_j,
 \qquad \eta^i{}_j~:=r^i{}_j-\fr1{q^2}v^iv_j,
 \qquad \eta^{ij}~:=r^{ij}-\fr1{q^2}v^iv^j ,
 \ee
which obey the identities
\be
\eta^n{}_j=a^{nm}\eta_{mj}, \quad  \eta^{ij}= a^{in}\eta_n^j,
\ee
\ses
\be
  \eta_{ni}y^i=0,
\ee
\ses
\be
 \eta_{ij}b^j=-(1-c^2)\fr1{q^2}z_i,
\qquad
 \eta_{ij}z^j=(1-c^2)\fr{S^2}{q^2}z_i,
 \ee
 \ses
and

\be
\D{\lf(\fr1qv_k\rg)}{y^j}=\fr1q\eta_{kj}, \qquad
\D{z_i}{y^k}=b\eta_{ik}
+\fr1{q^2}v_kz_i+\fr1b(b_kz_i-z_kb_i).
\ee

Using the generating metric function $V(x,w)$ defined from the representation
\be
K=bV(x,w),
\ee
we obtain
\be
\D K{y^i}=b_iV+\fr1{bq}z_iV', \qquad
\Dd{ K}{y^i}{y^j}=
\fr1q\eta_{ij}
V'
+\fr1{b^3q^2}z_iz_jV'',
\ee
where
\be
V'= \D Vw, \qquad V''= \DD Vw.
\ee

Taking into account the explicit derivatives of the function  $V$:
 \be
VV'=w\fr {K^2}{B},  \qquad    VV''= \fr{b^2}B  \fr {K^2}{B}
 \ee
(use (A.19)-(A.23)),
we find
the representations
\be
y_i=\Bigl(Bb_i+z_i\Bigr)\fr{K^2}{bB},
\ee
\ses
\be
g_{ij}=
\fr{K^2}B\eta_{ij}
+\fr{K^2}{b^2}b_ib_j
+
\fr{K^2}{b^2B}(b_iz_j+b_jz_i)
+\fr{B-gbq}{b^2q^2}\fr{K^2}{B^2}z_iz_j,
\ee
\ses
\ses
\be
h_{ij}=
\fr{K^2}B
\lf(\eta_{ij}
+\fr{1}{Bq^2}z_iz_j
\rg),
\ee
which entail
\be
y_i=\Bigl(v_i+(b+gq)b_i\Bigr)\fr{K^2}B,
\ee
\ses
\be
g_{ij}=
\biggl[a_{ij}
+\fr g{B}\Bigl (q(b+gq)b_ib_j+q(b_iv_j+b_jv_i)-b\fr{v_iv_j}q\Bigr)\biggr]\fr{K^2}B,
\ee
and
\be
g^{ij}=
\biggl[a^{ij}
+\fr g{B}\Bigl(-bqb^ib^j-q(b^iv^j+b^jv^i)+(b+gc^2q)\fr{v^iv^j}{\nu}\Bigr)
\biggr]\fr B{K^2}.
\ee
The determinant of the metric tensor  is everywhere positive:
\be
\det(g_{ij})=\fr {\nu}q\biggl(\fr{K^2}B\biggr)^N\det(a_{ij})>0
\ee
with
the function $\nu$ given by (A.29).

In terms of the set $\{b_i,u_i=a_{ij}y^j\}$, we obtain the alternative representations
\be
y_i=( u_i+gqb_i ) \fr{K^2}B,
\ee
\ses
\be
g_{ij}=
\biggl[a_{ij}
+\fr g{B}\Bigl((gq^2-\fr{bS^2}q)b_ib_j-\fr bqu_iu_j+
\fr{ S^2}q(b_iu_j+b_ju_i)\Bigr)\biggr]\fr{K^2}B,
\ee
and
\be
g^{ij}=
\biggl[a^{ij}+\fr g{\nu}(bb^ib^j-b^iy^j-b^jy^i)+\fr g{B\nu}(b+gc^2q)y^iy^j
\biggr]\fr B{K^2},
\ee
\ses
together with
\be
h_{ij}=
\biggl[a_{ij}
+\fr 1{qB}\Bigl(
-gbS^2b_ib_j-(q+gb)u_iu_j+gb^2(b_iu_j+b_ju_i)
\Bigr)\biggr]\fr{K^2}B,
\ee
\ses
which entails
\be
h_{ij}b^j=-\fr{\nu z_i}{Bq}\fr{K^2}B, \qquad
g^{ij}z_iz_j=\fr q{\nu}(c^2S^2-b^2)\fr{B^2}{K^2}
\ee
\ses
and
\be
g^{ij}b_j=\fr1{K^2}\Bigl[S^2b^i-\fr g{\nu}(c^2S^2-b^2)v^i\Bigr].
\ee

\ses

Given any vector $t_j$, we have
$$
g^{ij}t_j=
\Biggl[Ba^{ij}t_j+\fr g{\nu}\biggl(B(bb^ib^j-b^iy^j-b^jy^i)+(b+gc^2q)y^iy^j
\biggr)t_j
\Biggr]\fr 1{K^2},
$$
\ses
or
\be
g^{ij}t_j=
\Biggl[Ba^{ij}t_j
-gq(yt)b^i+\fr g{\nu}
\Bigl(-B(bt)+(b+gc^2q)(yt)\Bigr)v^i
\Biggr]\fr 1{K^2}.
\ee

Also,
\be
b+gc^2q=\fr1b(B-q\nu), \qquad
\D B{y^k}=\fr{2B}{K^2}y_k+\fr gq z_k,
\ee
\ses
\be
v_k=\fr{q^2}{K^2}y_k+\fr{b+gq}B z_k,
\ee
\ses
\be
g^{ij}u_j=\fr1{K^2}
\Biggl(
By^i-gq\Bigr[S^2b^i-\fr g{\nu}(c^2S^2-b^2)v^i \Bigr]
\Biggr),
\ee
\ses
and
\be
g^{ij}v_j=\fr1{K^2}
\Biggl
(By^i-(b+gq)\Bigr[S^2b^i-\fr g{\nu}(c^2S^2-b^2)v^i\Bigr]
\Biggr).
\ee

Using the function $X$ given by
\be
\fr1X=    N+\fr{(1-c^2)B}{q\nu},
\ee
we can evaluate the Cartan tensor
\be
 A_{ijk}~ := \fr K2\D{g_{ij}}{y^k}
 \ee
and the contraction
\be
 A_i~:=g^{jk}A_{ijk} =  K\D{\ln\bigl(\sqrt{\det(g_{mn})}\bigr)}{y^i}.
\ee
From (A.55) it follows that
\be
A_i=\fr {Kg}{2qB}
\fr1X
(q^2b_i- bv_i),
\ee
which entails
\be
A^iA_i=     \fr{g^2}4
\fr1{X^2}\lf(N+1-\fr1X\rg).
\ee
Also,
we find
\be
A_{ijk}= X
 \Biggl[
A_ih_{jk}  +A_jh_{ik}  +A_kh_{ij}
-\lf(N+1-\fr1X\rg)
\fr1{A_hA^h}A_iA_jA_k
\Biggr].
\ee

The equalities
\be
A^i=
\fr {g}{2XK\nu}
\Bigl[ Bb^i
-   (b+gqc^2)y^i
\Bigr],
\ee
\ses
\be
b_iA^i=   \fr {g}{2XK\nu}  (c^2S^2-b^2),
\ee
\ses
and
\be
v_i=\fr {q^2}{K^2}y_i  -   q(b+gq) \fr{2X}{Kg} A_i
\ee
are valid.

\ses

If we use the vector
\be
e_i~:=-b_i+\fr{bv_i}{q^2},
\ee
so that
$$
y^ie_i=0,
$$
we can readily convert the  representation (A.72) to the form
\be
\fr {B^2}{K^3}A_{ijk}=-\fr12gq
(
e_k\eta_{ij}+e_i\eta_{kj}+e_j\eta_{ik}
)
-\fr {gq^3}{B}e_ie_je_k.
\ee

In various processes of evaluations, it is useful to apply the formulas
\be
Kb_n= bl_n+ \fr{2q}{g}X A_n,
 \qquad
K\fr1q v_n=ql_n  -   \fr{B-q^2}b \fr{2}{g} XA_n,
\ee
\ses
and
\be
\D{\lf( \fr Kq\rg)}{y^n}=\fr2{gbq^2}(B-q^2) XA_n,
\ee
\ses
together with
\be
\D B{y^k}=\fr{2B}{K^2}y_k      -      \fr{2B}K XA_k,
\qquad
\D{\lf(\fr{q^2}B\rg)} {y^k}=  -      \fr{2q(2b+gq)}{gKB} XA_k.
\ee

The formula (A.67) can also be represented in the form
\be
\fr1X=N+ \fr{1-c^2}w \fr{1+gw+w^2}{w+(1-c^2)g}.
\ee

Simple straightforward calculation yields
\be
\D{(XA_k)}{y^n}  =  -      \fr {1}{K}l_k\, XA_n
-\fr{g}{2K w}  h_{kn}
+ \fr{2}{gKq}(b+gq)   X^2A_kA_n.
\ee

\ses

We use
 the Riemannian covariant derivative
\be
\nabla_ib_j~:=\partial_ib_j-b_ka^k{}_{ij},
\ee
where
\be
a^k{}_{ij}~:=\fr12a^{kn}(\prtl_ja_{ni}+\prtl_ia_{nj}-\prtl_na_{ji})
\ee
are the
Christoffel symbols given rise to by the associated Riemannian metric ${\cal S}$.

\ses\ses

\setcounter{equation}{0}

{
\nin \bf Appendix B:   Indefinite
 ${\mathbf\cF\cR^{SR}_{g;c} } $-space
 }

\ses\ses

The positive--definite ${\mathbf\cF\cR^{PD}_{g;c} }  $--space described
possesses
the indefinite (relativistic) version, to be denoted as  the ${\mathbf\cF\cR^{SR}_{g;c} }  $--space
(with the upperscripts  ``SR'' meaning ``special--relativistic'').
The underlined space $\cR_N=\{M, a_{mn}\}$ is now taken to be {\it pseudo-Riemannian},
such that   the input metric tensor $\{a_{mn}(x)\}$  is to be
 pseudo--Riemannian with
the {\it time--space signature:}
\be
{\rm sign}(a_{mn})=(+ - -\dots).
\ee
The definition range
$-2<g(x)<2 $ and the representation $h=\sqrt{1-(1/4)g^2}$
applicable in the positive-definite case
(see (A.11) and (A.12))
transform now according to
$$
 -\infty < g(x) < \infty, \qquad
h(x)=\sqrt{1+\fr14(g(x))^2}, \qquad
G(x)=\fr{g(x)}{h(x)}
 $$
(such a  phenomenon was explained in  [3]).
The {\it  pseudo--Finsleroid-regular characteristic
quadratic form}
\be
B(x,y) :=b^2-gqb-q^2
\equiv (b+g_+q)(b+g_-q)
\ee
is now of the positive discriminant
\be
D_{\{B\}}=4h^2>0
\ee
(compare these formulas with (A.13) and (A.14)).

In terms of these concepts, we propose

\ses

 { DEFINITION}. The scalar function $F(x,y)$ given by the formula
\be
F(x,y)~:=\sqrt{|B(x,y)|}\,J(x,y)
\equiv
|b+g_-q|^{G_+/2}|b+g_+q|^{-G_-/2},
\ee
where
\be
J(x,y)=
\lf|
\fr{b+g_-q}{b+g_+q}
\rg|^{-G/4},
\ee
is called
the {\it  pseudo-Finsleroid-regular  metric function}.
It is convenient to use the quantities
\be
g_+=-\fr12g+h, \qquad g_-=-\fr12g-h,
\ee
\medskip
\be
G_+=\fr{g_+}h\equiv -\fr12G+1, \qquad G_-=\fr{g_-}h\equiv -\fr12G-1.
\ee

\ses
\ses

Again,
the zero--vector $y=0$ is excluded from consideration:
$
y\ne 0.
$
The positive  (not absolute) homogeneity  holds:
\be
F(x,\la y)=\la F(x,y), \qquad \la>0, ~ \forall x, ~ \forall y.
\ee
The function
 \be
 L(x,y) =q-\fr g2b
\ee
is now to be used instead of (A.23), so that (A.24) changes to read
\be
L^2-h^2b^2=B.
\ee

{

Similarly to (A.25), we introduce

\ses

 { DEFINITION}.  The arisen  space
\be
\cF\cR^{SR}_{g;c} :=\{\cR_{N};\,b_i(x);\,g(x);\,F(x,y)\}
\ee
is called the
 {\it pseudo--Finsleroid-regular space}.

\ses

\ses

 { DEFINITION}. The space $\cR_N=(M,\,{\cal S})$ entering the above definition (B.11)
 is called the {\it associated  pseudo--Riemannian space}.

\ses

 { DEFINITION}. The scalar $g(x)$ is called
the {\it  pseudo--Finsleroid charge}.
The 1-form $b$ is called the  {\it  pseudo--Finsleroid--axis}  1-{\it form}.

\ses

The equality
\be
a_{ij}(x)=g_{ij}(x,y)\bigl|_{g=0}\bigr.
\ee
(cf. (A.26)) is applicable to the pseudo--Finsleroid case.

One can  observe the phenomenon that  the representations of the components
$y_i,\, g_{ij},\,g^{ij},\,A_i,\,A_{ijk}$
in
the $\cF\cR^{SR}_{g;c} $--space
are directly obtainable from the positive--definite case representations (written in the preceding Appendix A)
through the formal change:
\be
g ~ \stackrel{PD ~ \to ~SR}{\Longrightarrow} ~ ig
\ee
and
\be
q ~ \stackrel{PD ~ \to ~SR}{\Longrightarrow} ~ iq,
\ee
where $i$ stands for the imaginary unity.
Therefore, we may apply the rules
\be
\fr gq ~ \stackrel{PD ~ \to ~SR}{\Longrightarrow} ~ \fr gq,   \qquad
 gq ~ \stackrel{PD ~ \to ~SR}{\Longrightarrow} ~  - gq.
\ee
It is the useful exercise to verify that if we apply these rules
 to the expression (B.4) of the relativistic function $F$, we obtain the positive--definite
case function $K$ defined by (A.19).


\ses

\ses

\def\bibit[#1]#2\par{\rm\noindent\parskip1pt
                     \parbox[t]{.05\textwidth}{\mbox{}\hfill[#1]}\hfill
                     \parbox[t]{.925\textwidth}{\baselineskip11pt#2}\par}

\nin
{  REFERENCES}

\ses

\bibit[1] G.S. Asanov:  Finsleroid--Finsler  space with Berwald and  Landsberg conditions,
 {\it  arXiv:math.DG}/0603472 (2006).

\ses

\bibit[2] G.S. Asanov:  Finsleroid--Finsler  space and spray   coefficients, {\it  arXiv:math.DG}/0604526 (2006).

\ses

 \bibit[3] G.S. Asanov:  Finsleroid--Finsler  spaces of positive--definite and  relativistic types.
 \it Rep. Math. Phys. \bf 58 \rm(2006), 275--300.

\ses

\bibit[4] G. S. Asanov:  Finsleroid--Finsler space and geodesic spray    coefficients,
{\it Publ.  Math. Debrecen } {\bf 71/3-4} (2007), 397-412.

\ses

 \bibit[5] D.~Bao, S.S. Chern, and Z. Shen: {\it  An
Introduction to Riemann-Finsler Geometry,}  Springer, N.Y., Berlin 2000.

\ses

\bibit[6] H. Rund: \it The Differential Geometry of Finsler
 Spaces, \rm Springer, Berlin 1959.

\end{document}